\documentstyle[amscd,amssymb,verbatim,epsf]{amsart}
\begin{document}

\theoremstyle{plain}
\newtheorem{Thm}{Theorem}
\newtheorem{Cor}{Corollary}
\newtheorem{Con}{Conjecture}
\newtheorem{Main}{Main Theorem}

\newtheorem{Lem}{Lemma}
\newtheorem{Prop}{Proposition}

\theoremstyle{definition}
\newtheorem{Def}{Definition}
\newtheorem{Note}{Note}

\theoremstyle{remark}
\newtheorem{notation}{Notation}
\renewcommand{\thenotation}{}

\errorcontextlines=0
\numberwithin{equation}{section}
\renewcommand{\rm}{\normalshape}%

\title[Strongly Cesaro type quasi-Cauchy sequences]%
   {Strongly Cesaro type quasi-Cauchy sequences}
%\author{H\"usey\.{I}n \c{C}akall\i \\ Maltepe University, Istanbul, Turkey}
%%, Department of Mathematics, Maltepe University,TR 34857, Maltepe, Istanbul, Turkey Phone:(+90216)6261050 ext:1206, \;  fax:(+90216)6261113
\author{H\"usey\.{I}n \c{C}akall\i\ \\ Maltepe University, Istanbul, Turkey Phone:(+90216)6261050/2248, fax:(+90216)6261113}

\address{H\"usey\.{I}n \c{C}akall\i \\
          Maltepe University, Department of Mathematics, Marmara E\u{g}\.{I}t\.{I}m K\"oy\"u, TR 34857, Maltepe, \.{I}stanbul-Turkey \; \; \; \; \; Phone:(+90216)6261050 ext:2248, \;  fax:(+90216)6261113}

\email{hcakalli@@maltepe.edu.tr; hcakalli@@gmail.com}

\keywords{Strongly Cesaro summable sequences, $N_{\theta}$-convergence, statistical convergence, quasi-Cauchy sequences, continuity}
\subjclass[2000]{Primary: 40A05 Secondaries: 40D25, 26A15}
\date{\today}

%$\boldsymbol{\alpha}=(\alpha_{n})$

\maketitle

\begin{abstract}
In this paper we call a real-valued function $N_{\theta}$-ward continuous if it preserves $N_{\theta}$-quasi-Cauchy sequences where a sequence $\boldsymbol{\alpha}=(\alpha _{k})$ is defined to be $N_{\theta}$-quasi-Cauchy when the sequence $\Delta \boldsymbol{\alpha}$ is in $N^{0}_{\theta}$. We prove not only inclusion and compactness type theorems, but also continuity type theorems.

\end{abstract}

\maketitle

\section{Introduction}

\normalfont{}
The concept of continuity and any concept involving continuity play a very important role not only in pure mathematics but also in other branches of sciences involving mathematics especially in computer science, information theory, biological science.

A real function $f$ is continuous if and only if, for each point $\alpha _{0}$ in the domain, $\lim_{n\rightarrow\infty}f(\alpha _{n})=f(\alpha _{0})$ whenever $\lim_{n\rightarrow\infty}\alpha _{n}=\alpha _{0}$. This is equivalent to the statement that $(f(\alpha _{n}))$ is a convergent sequence whenever $(\alpha _{n})$ is. This is also equivalent to the statement that $(f(\alpha _{n}))$ is a Cauchy sequence whenever $(\alpha _{n})$ is provided that the domain of the function is either whole $\Re$ or a bounded and closed subset of $\Re$ where $\Re$ is the set of real numbers. These well known results for continuity for real functions in terms of sequences might suggest to us giving a new type of continuity, namely, $N_{\theta}$-ward continuity.

The purpose of this paper is to introduce a concept of $N_{\theta}$-ward continuity of a real function and a concept of $N_{\theta}$-ward compactness of a subset of $\Re$ which cannot be given by means of any $G$ and prove some theorems.

\maketitle

\section{Preliminaries}

\normalfont{}

We will use boldface letters $\boldsymbol{\alpha}$, $\bf{x}$, $\bf{y}$, $\bf{z}$, ... for sequences $\boldsymbol{\alpha}=(\alpha_{k})$, $\textbf{p}=(p_{n})$, $\textbf{x}=(x_{n})$, $\textbf{y}=(y_{n})$, $\textbf{z}=(z_{n})$, ... of points in $\Re$ for the sake of abbreviation. $s$ and $c$ will denote the set of all sequences, and the set of convergent sequences of points in $\Re$.

 A subset of $\Re$ is compact if and only if it is closed and bounded. A subset $A$ of  $\Re$ is bounded if $|a|\leq{M}$ for all $a \in{A}$  where  $M$ is a positive  real constant number. This is equivalent to the statement that any sequence of points in $A$ has a Cauchy subsequence. The concept of a Cauchy sequence involves far more than that the distance between successive terms is tending to zero. Nevertheless, sequences which satisfy this weaker property are interesting in their own right.  A sequence $(\alpha _{n})$ of points in $\Re$ is quasi-Cauchy if $(\Delta \alpha _{n})$ is a null sequence where $\Delta \alpha _{n}=\alpha _{n+1}-\alpha _{n}$. These sequences were named as quasi-Cauchy by Burton and Coleman \cite{BurtonColeman}, while they  were called as forward sequences in \cite{CakalliForwardcontinuity},(see also \cite{KizmazOncertainsequencespaces}, and \cite{EmrahEvrenKaraandBasarirOncompactoperatorsandsomeEulerB(m)differencesequencespaces}).

It is known that a sequence $(\alpha _{n})$ of points in $\Re$ is slowly oscillating if $$
\lim_{\lambda \rightarrow 1^{+}}\overline{\lim}_{n}\max _{n+1\leq
k\leq [\lambda n]} |
  \alpha _{k}  -\alpha _{n} | =0
$$ where $[\lambda n]$ denotes the integer part of $\lambda n$. This is equivalent to the following if $\alpha _{m}-\alpha _{n}\rightarrow 0$ whenever $1\leq \frac{m}{n}\rightarrow 1$ as, $m,n\rightarrow \infty$.
Using $\varepsilon>0$ s and $\delta$ s this is also equivalent to the case when for any given $\varepsilon>0$, there exist $\delta=\delta (\varepsilon) >0$ and a positive integer $N=N(\varepsilon)$ such that $|\alpha _{m}-\alpha _{n}|<\varepsilon$ if $n\geq N(\varepsilon)$ and $n\leq m \leq (1+\delta)n$ (see \cite{FDikMDikandCanak}).
Any Cauchy sequence is slowly oscillating, and any slowly oscillating sequence is quasi-Cauchy. There are quasi-Cauchy sequences which are not Cauchy. For example, the sequence  $(\sqrt{n})$ is quasi-Cauchy, but Cauchy. Any subsequence of a Cauchy sequence is Cauchy. The analogous property fails for quasi-Cauchy sequences, and fails for slowly oscillating sequences as well. A counterexample for the case, quasi-Cauchy, is again the sequence $(a_n)=(\sqrt{n})$ with the subsequence $(a_{n^{2}})=(n)$. A counterexample for the case slowly oscillating is the sequence $(log_{10} n)$ with the subsequence (n). Furthermore we give more examples without neglecting (see \cite{Vallin}): the sequences $(\sum^{\infty}_{k=1}\frac{1}{n})$, $(ln\;n)$, $(ln\;(ln\;n))$, $(ln\;( ln\;( ln n)))$ and combinations like that are all slowly oscillating, but Cauchy. The bounded sequence $(cos (6 log(n + 1)))$ is slowly oscillating, but Cauchy. The sequences $(cos(\pi \sqrt{n}))$ and $(\sum^{k=n}_{k=1}(\frac{1}{k})(\sum^{j=k}_{j=1}\frac{1}{j}))$ are quasi-Cauchy, but slowly oscillating (see also  \cite{CakalliSlowlyoscillatingcontinuity}, \cite{DikandCanak},  and \cite{CakalliNewkindsofcontinuities}).

By a method of sequential convergence, or briefly a method, we mean a linear function $G$ defined on a subspace of $s$, denoted by $c_{G}$, into $\Re$. A sequence $\textbf{x}=(x_{n})$ is said to be $G$-convergent to $\ell$ if $\textbf{x}\in c_{G}$ and $G(\textbf{x})=\ell$ (see \cite{CakalliSequentialdefinitionsofcompactness}). In particular, $\lim$ denotes the limit function $\lim \textbf{x}=\lim_{n}x_{n}$ on the subspace $c$. A method $G$ is called regular if $c\subset{c_{G}}$, i.e. every convergent sequence $\textbf{x}=(x_{n})$ is $G$-convergent with $G(\textbf{x})=\lim \textbf{x}$.

Now we discuss some special classes of methods of sequential convergence that have been studied in the literature for real or complex number sequences. Firstly, for real and complex number sequences, we note that the most important transformation class is the class of matrix methods. Consider an infinite matrix $\verb"A"=(a_{nk})^{\infty}_{n,k=1}$ of real numbers. Then, for any sequence $\textbf{x}=(x_{n})$ the sequence $\verb"A"\textbf{x}$ is defined as

\[
\verb"A"\textbf{x}=(\sum^{\infty}_{k=1}a_{nk}x_{k})_{n}
\]
provided that each of the series converges. A sequence $\textbf{x}$ is called $\verb"A"$-convergent (or $\verb"A"$-summable) to $\ell$ if $\verb"A"\textbf{x}$ exists and is convergent with

\[
\lim \verb"A"\textbf{x}=\lim_{n\rightarrow\infty}\sum^{\infty}_{k=1}a_{nk}x_{k}=\ell.
\]
Then $\ell$ is called the $\verb"A"$-limit of $\textbf{x}$. We have thus defined a method of sequential convergence, i.e. $G(\textbf{x})=\lim \verb"A"\textbf{x}$, called a matrix method or a summability matrix. A subset $A$ of $\Re$ is $G$-sequentially compact if any sequence $\textbf{x}$ of points in $A$ has a subsequence $\textbf{z}$ such that $G(\textbf{z})\in {A}$.

The Hahn-Banach theorem can be used to define methods which are not generated by a regular summability matrix. Banach used this theorem to show that the limit functional can be extended from the convergent sequences to the bounded sequences while preserving linearity, positivity and translation invariance; these extensions have come to be known as Banach limits. If a bounded sequence is assigned the same value  $\ell$  by each Banach limit, the sequence is said to be almost convergent to $\ell$. It is well known that a sequence $\boldsymbol{\alpha}=(\alpha_{n})$ is almost convergent to $\ell$ if and only if
\[
\lim_{n\rightarrow\infty}\frac{1}{n}\sum^{n}_{k=1}\alpha_{k+j}=\ell,
\]
uniformly in $j$. This defines a method of sequential convergence, i.e. \\$G(\boldsymbol{\alpha}):=almost\;limit\;of\;\textbf{x}$.

%We define a concept of almost quasi-Cauchy of a sequence
We call a sequence $\boldsymbol{\alpha}=(\alpha_{n})$ almost quasi-Cauchy if
\[
\lim_{n\rightarrow\infty}\frac{1}{n}\sum^{n}_{k=1} (\alpha_{k+j}-\alpha_{k+j+1})=0,
\]
uniformly in $j$. Any almost convergent sequence is almost quasi-Cauchy, but the converse is not always true.

The idea of statistical convergence was formerly given under the name "almost convergence" by Zygmund in the first edition of his celebrated monograph published in Warsaw in 1935 \cite{ZygmundTrigonometricseries}. The concept was formally introduced by Fast \cite{Fast} and later was reintroduced by Schoenberg \cite{SchoenbergTheintegrabilityofcertainfunctions}, and also independently by Buck \cite{BuckGeneralizedasymptoticdensity}. Although statistical convergence was introduced over nearly the last seventy years, it has become an active area of research for twenty years. This concept has been applied in various areas such as number theory \cite{ErdosSurlesdensitesde}, measure theory \cite{MillerAmeasuretheoresubsequencecharacterizationofstatisticalconvergence}, trigonometric series \cite{ZygmundTrigonometricseries}, summability theory \cite{FreedmanandSemberDensitiesandsummability}, locally convex spaces \cite{MaddoxStatisticalconvergenceinlocallyconvex}, in the study of strong integral summability \cite{ConnorandSwardsonStrongintegralsummabilityandstonecompactification}, turnpike theory \cite{MakarovLevinRubinovMathematicalEconomicTheory}, \cite{MckenzieTurnpiketheory}, \cite{PehlivanandMamedovStatisticalclusterpointsandturnpike}, Banach spaces \cite{ConnorGanichevandKadetsAcharacterizationofBanachspaceswithseparabledualsviaweakstatisticalconvergence}, and metrizable topological groups \cite{Cakallilacunarystatisticalconvergenceintopgroups}, and topological spaces \cite{MaioKocinac}, \cite{CakalliandKhan}. It should be also mentioned that the notion of statistical convergence has been considered, in other contexts, by several people like R.A. Bernstein, Z. Frolik, etc. The concept of statistical convergence is a generalization of the usual notion of convergence that, for real-valued sequences, parallels the usual theory of convergence. A sequence $(\alpha _{k})$ of points in $\Re$ is called statistically convergent to an element $\ell$ of $\Re$  if for each
$\varepsilon$
\[
\lim_{n\rightarrow\infty}\frac{1}{n}|\{k\leq n: |\alpha _{k}-\ell|\geq{\varepsilon}\}|=0,
\] and this is denoted by $st-\lim_{k\rightarrow\infty}\alpha _{k}=\ell$ \cite{Fridy} (see also \cite{CakalliAstudyonstatisticalconvergence}).

A sequence $\boldsymbol{\alpha}=(\alpha_{n})$ of real numbers is called Abel
convergent (or Abel summable)  to $\ell$ if
the series $\Sigma_{k=0}^{\infty}\alpha_{k}x^{k}$ is convergent for $0\leq x<1$ and
\[\lim_{x \to 1^{-}}(1-x) \sum_{k=0}^{\infty}\alpha_{k}x^{k}=\ell.\]
In this case we write  $Abel-\lim \alpha_{n}=\ell$.
The set of Abel convergent sequences will be denoted by $\textbf{A}$.
Abel proved that if $\lim_{n \to \infty}\alpha_{n}=\ell$, then $Abel-\lim \alpha_{n}=\ell$, i.e. every convergent sequence is Abel convergent to the same limit (\cite{Abel}, see also \cite{CVStanojevicandVBStanojevic}, and \cite{FridyandKhan}). As it is known that the converse is not always true in general, as we see that the sequence $((-1)^{n})$ is Abel convergent to $0$ but convergent in the ordinary sense (see \cite{CakalliandAlbayrakAbelContinuity}).

Now we recall the definitions of ward compactness, and slowly oscillating compactness.

\begin{Def}
(\cite{CakalliForwardcontinuity}) A subset $A$ of  $\Re$ is called ward compact if whenever $(\alpha _{n})$ is a sequence of points in $A$ there is a quasi-Cauchy subsequence $\textbf{z}=(z_{k})=(\alpha_{n_{k}})$ of $(\alpha _{n})$.

\end{Def}

\begin{Def}
(\cite{CakalliSlowlyoscillatingcontinuity}, \cite{CakalliNewkindsofcontinuities}) A subset $A$ of  $\Re$ is called slowly oscillating compact if whenever $(\alpha _{n})$ is a sequence of points in $A$ there is a slowly oscillating subsequence $\textbf{z}=(z_{k})=(\alpha_{n_{k}})$ of $(\alpha _{n})$.\\

\end{Def}

In an unpublished work, \c{C}akall\i\; called a sequence $(\alpha _{k})$ of points in $\Re$ statistically quasi-Cauchy if $$st-\lim_{k\rightarrow\infty}\Delta \alpha _{k}=0.$$

Any quasi-Cauchy sequence is statistically quasi-Cauchy, but the converse is not always true. Any statistically convergent sequence is statistically quasi-Cauchy. There are statistically quasi-Cauchy sequences which are not statistically convergent.

Now we recollect the following two definitions of an unpublished work of \c{C}akall\i\;.

\begin{Def}

 A subset $A$ of  $\Re$ is called statistically ward compact if whenever $(\alpha _{n})$ is a sequence of points in $A$ there is a statistical quasi-Cauchy subsequence $\textbf{z}=(z_{k})=(\alpha _{n_{k}})$ of $(\alpha _{n})$.

\end{Def}

\begin{Def}

A function defined on a subset $A$ of $\Re$ is called statistically ward continuous if it preserves statistically quasi-Cauchy sequences, i.e. $(f(\alpha_{n}))$ is a statistically quasi-Cauchy sequence whenever $(\alpha_{n})$ is.
\end{Def}

Now we give the following definition.

\begin{Def}
 A subset $A$ of  $\Re$ is called almost ward compact if whenever $(\alpha _{n})$ is a sequence of points in $A$ there is a almost quasi-Cauchy subsequence $\textbf{z}=(z_{k})=(\alpha_{n_{k}})$ of $(\alpha _{n})$.

\end{Def}

A sequence $\boldsymbol{\alpha}=(\alpha_{k})$ is called strongly Cesaro summable to a real number $\ell$ if $\lim_{n\rightarrow\infty}\frac{1}{n}\sum^{n}_{i=1}|\alpha_{i}-\ell|=0 $. This is denoted by $|C_{1}|-lim \alpha_{k}=\ell$. The set of all strongly Cesaro summable sequences is denoted by $|\sigma_{1}|$. A sequence $\boldsymbol{\alpha}=(\alpha_{k})$ is called strongly Cesaro quasi-Cauchy if $\lim_{n\rightarrow\infty}\frac{1}{n}\sum^{n}_{i=1}|\alpha_{i}-\alpha_{i+1}|=0 $. The set of all strongly Cesaro quasi-Cauchy sequences is denoted by $|\Delta \sigma_{1}|$.
%$|\Delta \sigma_{1}|$.
%The well-known space $|\Delta \sigma_{1}|$ of strongly Cesaro summable sequences is defined as \[|\Delta \sigma_{1}|=\{\boldsymbol{\alpha}=(\alpha_{k}):there\; exists\; an\; \ell\; such\; that\; \lim_{n\rightarrow\infty}\frac{1}{n}\sum^{n}_{i=1}|\alpha_{i}-\ell|=0 \}.\]

By a lacunary sequence $\theta=(k_{r})$, we mean an increasing sequence $\theta=(k_{r})$ of positive integers such that $k_{0}=0$ and $h_{r}:k_{r}-k_{r-1}\rightarrow\infty$. The intervals determined by $\theta$ will be denoted by $I_{r}=(k_{r-1}, k_{r}]$, and the ratio $\frac{k_{r}}{k_{r-1}}$ will be abbreviated by $q_{r}$.

Sums of the form $\sum^{k_{r}}_{k_{r-1}+1}|\alpha_{k}|$ frequently occur, and will often be written for convenience as $\sum^{}_{k\in{I_{r}}}|\alpha_{k}|$.

The notion of $N_\theta$ convergence was introduced, and studied by Freedman, Sember, and M. Raphael in \cite{FreedmanandSemberandRaphaelSomecesarotypesummabilityspaces}.
%Fridy and Orhan in \cite{FridyandOrhanlacunarystatisconvergence} and \cite{FridyandOrhanlacunarystatisticalsummability}
%(see also \cite{FreedmanandSemberandRaphaelSomecesarotypesummabilityspaces}) |\{k\in I_{r}: |\alpha_{k}-\ell| \geq\varepsilon\}|=0
 %and \cite{Cakallilacunarystatisticalconvergenceintopgroups}). \sum^{\infty}_{k\in{I_{r}}}|\alpha_{k}-\ell|=0
A sequence $(\alpha_{k})$ of points in $\Re$ is called $N_\theta$-convergent to an element $\ell$ of \textbf{R} if
\[
\lim_{r\rightarrow\infty}\frac{1}{h_{r}}\sum^{}_{k\in{I_{r}}}|\alpha_{k}-\ell|=0,
\]
and it is denoted by $N_{\theta}-lim\; \alpha_{k}=\ell$.
%for every positive real number $\varepsilon$.
%where $I_{r}=(k_{r-1},k_{r}]$ and $k_{0}=0$, $h_{r}:k_{r}-k_{r-1}\rightarrow \infty$ as $r\rightarrow\infty$ and $\theta=(k_{r})$ is an increasing sequence of positive integers. This is denoted by $N_\theta-\lim_{k\rightarrow\infty}\alpha _{k}=\ell$.
%$\liminf_{r}\frac{k_{r}}{k_{r-1}}>1$
Throughout the paper $c$, $S$, $N_{\theta}$ will denote the set of convergent sequences, the set of statistical convergent sequences, the set of $N_\theta$ convergent sequences of points in \textbf{R}, respectively.

Connor and Grosse-Erdmann \cite{ConnorandGrosse} gave sequential definitions of continuity for real functions calling it $G$-continuity instead of $A$-continuity and their results cover the earlier works related to $A$-continuity. In particular, $\lim$ denotes the limit function $\lim \textbf{x}=\lim_{n}x_{n}$ on the subspace $c$, and $st-lim$ denotes statistical limit function $st-lim \textbf{x}=st-lim x_{n}$ on the subspace, $S$, and $N_{\theta}-lim$ denotes $N_{\theta}$-limit function $N_{\theta}-lim \textbf{x}=N_{\theta}-lim x_{n}$ on the subspace, $N_{\theta}$. We see that $lim$, $st-lim$, and $Abel-lim$   are all regular methods without any restriction. $N_{\theta}-lim$  is regular under the condition $\lim inf_{r}\; q_{r}>1$.
 %which we assume throughout the paper.
 A function $f$ is called $G$-sequentially continuous at $u\in {\Re}$ provided that whenever a sequence $\textbf{x}$ of terms in $\Re$, then the sequence $f(\textbf{x})=(f(x_{n}))$ is $G$-convergent to $f(u)$.

Recently, \c{C}akall\i\; gave a sequential definition of compactness, which is a generalization of ordinary sequential compactness, as in the following: a subset $A$ of $\Re$ is $G$-sequentially compact if for any sequence $(\alpha _{k})$ of points in $A$ there exists a subsequence $\textbf z$ of the sequence such that $G(\textbf{z})\in{A}$. His idea enables to obtain new kinds of compactness via most of the non-matrix sequential convergence methods, for example Abel method, as well as all matrix sequential convergence methods.

\maketitle

\section{$N_{\theta}$-summable sequences}

We call a sequence $(\alpha _{k})$ of points in $\Re$, $N_{\theta}$-convergent to $\ell$ if $$N_{\theta}-\lim_{k\rightarrow\infty} \alpha _{k}=\ell.$$ Any convergent sequence is $N_{\theta}$-convergent, but the converse is not always true.

Throughout the paper $|\sigma_{1}|$ and $N_{\theta}$ will denote the set of strongly Cesaro convergent sequences and the set of $N_{\theta}$-convergent sequences of points in $\Re$, respectively. In this section, we will assume that $\lim inf_{r}\; q_{r}>1$.

Now we give the definition of $N_{\theta}$-sequentially compactness.

\begin{Def}
 A subset $A$ of  $\Re$ is called $N_{\theta}$-sequentially compact if whenever $(\alpha _{n})$ is a sequence of points in $A$ there is a $N_{\theta}$-convergent  subsequence $\textbf{z}=(z_{k})=(a_{n_{k}})$ of $(\alpha _{n})$ whose $N_{\theta}$-limit in $A$ .

\end{Def}

\noindent{\it \textbf{Lemma 1.}}
The sequential method $N_{\theta}$ is subsequential.

\noindent{\it \textbf{Theorem 2.}} A subset $A$ of $\Re$ is sequentially compact if and only if it is $N_{\theta}$-sequentially compact.

\begin{pf}

Let $A$ be any bounded closed subset of $\Re$ and  $(\alpha_{n})$ be any sequence of points in $A$.

As  $A$ is sequentially compact there is a convergent subsequence $(\alpha_{n_{k}})$ of $(\alpha_{n})$ whose limit in $A$. This subsequence is $N_{\theta}$-convergent since $N_{\theta}$-method is regular.

Thus (a) implies (b). The proof of the sufficiency follows from Lemma 1. Thus the proof of the theorem is completed.

\end{pf}

%\noindent{\it \textbf{Corollary 6.}} A closed subset of $\Re$ is $N_{\theta}$-sequentially compact if and only if it is $N_{\theta}$-sequentially compact.

%\noindent{\it \textbf{Corollary 7.}} A subset of $\Re$ is $N_{\theta}$-sequentially compact if and only if it is statistically ward compact.

A sequence $\alpha=(\alpha _{n})$ is $\delta$-quasi-Cauchy if $\lim_{k\rightarrow \infty} \Delta^{2} \alpha _{n}=0$ where $\Delta^{2}\alpha _{n}=a_{n+2}-2a_{n+1}+\alpha _{n}$ (\cite{CakalliDeltaquasiCauchysequences}). A subset $A$ of $\Re$ is called $\delta$-ward compact if whenever $\boldsymbol \alpha=(\alpha _{n})$ is a sequence of points in $A$ there is a subsequence $\textbf{z}=(z_{k})=(\alpha_{n_{k}})$ of $\boldsymbol \alpha$ with $\lim_{k\rightarrow \infty} \Delta^{2}z_{k}=0$. It follows from the above theorem that any $N_{\theta}$-sequentially compact subset of $\Re$ is $\delta$-ward compact.

We see that for any regular subsequential method $G$ defined on $\Re$, if a subset $A$ of $\Re$ is $G$-sequentially compact, then it is $N_{\theta}$-sequentially compact. But the converse is not always true.

\begin{Def}

 A subset $A$ of  $\Re$ is called strongly Cesaro compact if whenever $(\alpha _{n})$ is a sequence of points in $A$ there is a strongly Cesaro summable  subsequence $\textbf{z}=(z_{k})=(\alpha _{n_{k}})$ of $(\alpha _{n})$ to a strongly Cesaro limit that belongs to $A$.

\end{Def}

\begin{Def}

A function defined on a subset $A$ of $\Re$ is called strongly Cesaro continuous if it preserves strongly Cesaro summable sequences, i.e. $(f(\alpha_{n}))$ is a strongly Cesaro quasi-Cauchy sequence whenever $(\alpha_{n})$ is.
\end{Def}

In connection with $N_{\theta}$-convergent sequences and convergent sequences the problem arises to investigate the following types of  "continuity" of functions on $\Re$ under the condition that $\lim inf_{r}\; q_{r}>1$.

\begin{description}
\item[($ N_{\theta}$)] $(\alpha _{n}) \in { N_{\theta}} \Rightarrow (f(\alpha  _{n})) \in { N_{\theta}}$
\item[($ N_{\theta} c$)] $(\alpha  _{n}) \in { N_{\theta}} \Rightarrow (f(\alpha  _{n})) \in {c}$
\item[$(c)$] $(\alpha  _{n}) \in {c} \Rightarrow (f(\alpha  _{n})) \in {c}$
\item[$(c  N_{\theta})$] $(\alpha  _{n}) \in {c} \Rightarrow (f(\alpha  _{n})) \in { N_{\theta}}$
\item[($N_{\theta}$)] $(\alpha _{n}) \in {N_{\theta}} \Rightarrow (f(\alpha  _{n})) \in {N_{\theta}}$
\end{description}

We see that $( N_{\theta})$ is $N_{\theta}$-continuity of $f$, $(N_{\theta})$ is $N_{\theta}$-continuity of $f$, and $(c)$ is the ordinary continuity of $f$. It is easy to see that $(\delta  N_{\theta} c)$ implies $( N_{\theta})$, and $( N_{\theta} )$ does not imply $( N_{\theta} c)$;  and $( N_{\theta})$ implies $(c N_{\theta})$, and $(c  N_{\theta})$ does not imply $( N_{\theta})$; $( N_{\theta} c)$ implies $(c)$, and $(c)$ does not imply $( N_{\theta} c)$; and $(c)$ is equivalent to $(c  N_{\theta})$.

Now we give the implication $( N_{\theta})$ implies $(N_{\theta})$, i.e. any $N_{\theta}$-continuous function is $N_{\theta}$-continuous.

\noindent{\it \textbf{Lemma 3.}} A regular sequential method $G$ is subsequential if and only if $\overline {A}=\overline{A}^{G}$ for every subset $A$ of $\Re$.

\noindent{\it \textbf{Lemma 4.}} If a function $f$  is $G$-sequentially continuous on $\Re$, then $f(\overline{A}^{G})\subset {\overline {f(A)}^{G}}$ for any  subset $A$ of $\Re$.

\noindent{\it \textbf{Theorem 5.}} Let $\lim inf_{r}\; q_{r}>1$. A function $f$  is $N_{\theta}$-continuous on a subset $A$ of $\Re$ if and only if it is continuous on $A$.

\begin{pf}

The proof easily follows from Lemma 1, Lemma 3, and Lemma 4,  so is omitted.

\end{pf}

%We state the following straightforward result related to ordinary continuity.

%\noindent{\it \textbf{Corollary 9.}} If $f$ is $N_{\theta}$-continuous, then it is continuous in the ordinary sense.

%\noindent{\it \textbf{Corollary 10.}} If $f$ is $N_{\theta}$-continuous, then it is statistically continuous.

%Related to $G$-continuity we have much more in the following.

\noindent{\it \textbf{Corollary 6.}} Let $\lim inf_{r}\; q_{r}>1$. A real function is $N_{\theta}$-continuous if and only if it is $G$-continuous for any regular subsequential method $G$.

\maketitle

\section{$N_{\theta}$-quasi-Cauchy sequences}

We call a sequence $(\alpha _{k})$ of points in $\Re$ $N_{\theta}$-quasi-Cauchy if $$N_{\theta}-\lim_{k\rightarrow\infty}\Delta \alpha _{k}=0.$$ Any quasi-Cauchy sequence is $N_{\theta}$-quasi-Cauchy, but the converse is not always true. We see that a convergent sequence is $N_{\theta}$-quasi-Cauchy. There are $N_{\theta}$-quasi-Cauchy sequences which are not $N_{\theta}$-convergent. Throughout the paper $\Delta \sigma^{0}_{1}$ and $\Delta N^{0}_{\theta}$ will denote the set of strongly Cesaro quasi-Cauchy sequences and the set of $N_{\theta}$-quasi-Cauchy sequences of points in $\Re$, respectively.

Here we give some inclusion properties between the set of statistically quasi-Cauchy sequences and the set of $N_{\theta}$-quasi-Cauchy sequences.

\noindent{\it \textbf{Theorem 7.}} For any lacunary sequence $\theta$, $|\Delta \sigma^{0}_{1}| \subset {\Delta N^{0}_{\theta}}$ if and only if $$lim\; inf \;q_{r} > 1.$$
\begin{pf}

%The proof of Lemma 2.1 on page 509 in \cite{FreedmanandSemberandRaphaelSomecesarotypesummabilityspaces} cab be modified to the proof of this theorem, so is omitted.
Sufficiency. Let us first suppose that $lim\; inf \;q_{r} > 1$, $lim\; inf \;q_{r}=a$, say. Write $b=\frac{a-1}{2}$. Then there exists a positive integer $N$ such that $q_{r}\geq 1+b$ for $r\geq N$. Since  $h_{r}=k_{r}-k_{r-1}$, we have \\ $\frac{k_{r}}{h_{r}}\leq \frac{b+1}{b}$ and $\frac{k_{r-1}}{h_{r}}\leq \frac{1}{b}$ for $r\geq N$.\\ Let $(\alpha _{k})\in{|\Delta \sigma^{0}_{1}|}$. Write $t_{r}=\frac{1}{h_{r}}\sum^{}_{k\in{I_{r}}}|\Delta \alpha_{k}|$. Then \[t_{r}=\frac{k_r}{h_r}(\frac{1}{k_{r}}\sum^{k_r}_{i=1}|\Delta \alpha_{i}|)-\frac{k_{r-1}}{h_r}(\frac{1}{k_{r-1}}\sum^{k_{r-1}}_{i=1}|\Delta \alpha_{i}|).\] The sequences $(\frac{1}{k_{r}}\sum^{k_r}_{i=1}|\Delta \alpha_{i}|)$ and $(\frac{1}{k_{r-1}}\sum^{k_{r-1}}_{i=1}|\Delta \alpha_{i}|)$ both converge to $0$, and it follows that that $(\alpha _{k})\in{\Delta N^{0}_{\theta}}$.

Necessity. Now let us suppose that $lim\; inf \;q_{r} = 1$. Then we can choose a subsequence $(k_{r_{j}})$ of the lacunary sequence $\theta$ such that \\ $$\frac {k_{r_{j}}}{k_{r_{j}-1}}<1+\frac{1}{j}$$ and $$\frac{k_{r_{j}-1}}{k_{r_{j-1}}}>j$$ where $r_{j}\geq r_{j-1}+2$. Now define a sequence $(\alpha_{k})$ by $\alpha_{k}=1+\frac{1+(-1)^{k}}{2}$ if $k\in{I_{r_{j}}}$ for some $j=1,2,...,n,..., $ and $\alpha_{k}=0$ otherwise. It follows that $\alpha \notin {N^{0}_{\theta}}$. However the sequence $(\alpha_{k})$ defined in this way is strongly Cesaro-quasi-Cauchy.

\end{pf}
\noindent{\it \textbf{Theorem 8.}}
For any lacunary sequence $\theta$, $\Delta N^{0}_{\theta} \subset {|\Delta \sigma^{0}_{1}|}$ if and only if $$lim\; sup \;q_{r} < \infty.$$
\begin{pf}
%The proof of Lemma 2.2 on page 510 in \cite{FreedmanandSemberandRaphaelSomecesarotypesummabilityspaces} cab be modified to the proof of this theorem, so is omitted.
Sufficiency. If $lim\; sup \;q_{r} < \infty$, then we can find an $H>0$ such that $q_{r}<H$ for any positive integer $r$. Let $\boldsymbol{\alpha} \in \Delta N^{0}_{\theta}$ and $\varepsilon>0$. Then there exist an $n_{0} \in {N}$ and a positive real number $K$ such that $\tau_{i} < \varepsilon$ for $i \geq n_{0}$  and  $\tau_{i} \leq {K}$ for all $i \in {N}$. Then if $n$ is any positive integer with $k_{r-1}<n\leq {k_r}$, where $r>n_{0}$, we have   $\frac{1}{n}\sum^{n}_{i=1}|\Delta \alpha_{i}|\leq  \frac{1}{k_{r-1}} \sum^{k_{r}}_{i=1}|\Delta \alpha_{i}|=\frac{1}{k_{r-1}} (\sum^{}_{I_{1}}|\Delta \alpha_{i}|+\sum^{}_{I_{2}}|\Delta \alpha_{i}|+...+\sum^{}_{I_{r}}|\Delta \alpha_{i}| )$\\ $=\frac{k_{1}}{k_{r-1}}\tau_{1}+\frac{k_{2}-k_{1}}{k_{r-1}}\tau_{2}+...+\frac{k_{n_{0}}-k_{n_{0}-1}}{k_{r-1}}\tau_{n_{0}}+\frac{k_{n_{0}+1}-k_{n_{0}}}{k_{r-1}}\tau_{n_{0}+1} +...+\frac{k_{r}-k_{r-1}}{k_{r-1}}\tau_{r} $\\    $\leq K  \frac{k_{n_{0}}}{k_{r-1}}+\frac{k_{r}-k_{n_{0}}}{k_{r-1}}\varepsilon < \frac{k_{n_{0}}}{k_{r-1}}+H \varepsilon $.\\ Now it follows that $\frac{1}{n}\sum^{\infty}_{i=1}|\Delta \alpha_{i}|\leq  \frac{1}{k_{r-1}} \sum^{k_{r}}_{i=1}|\Delta \alpha_{i}|=0$. Hence $\boldsymbol{\alpha} \in {|\Delta \sigma^{0}_{1}|}$.

 Necessity. To prove that if  $\Delta N_{\theta} \subset {\Delta S }$, then $lim\; sup \;q_{r} < \infty$, suppose that $lim\; sup \;q_{r} = \infty$. Let $c$ be a fixed positive real number. Select a subsequence $(k_{r_{j}})$ of the lacunary sequence $\theta=(k_{r})$ such that $q_{r_{j}}>j$, $k_{r_{j}}>j+3$, and define a sequence $(\alpha_{k})$ by $\alpha_{k}=c+\frac{c+(-1)^{k}c}{2}$ if $k_{r_{j}-1}<k\leq 2k_{r_{j}-1}$ for some $j=1,2,...,$ and $\alpha _{k}=0$ otherwise. An easy calculation shows that $\boldsymbol{\alpha}\in{\Delta N^{0}_{\theta}}$. However $\boldsymbol{\alpha}\notin{|\Delta \sigma^{0}_{1}|}$. This completes the proof of necessity hance the proof of the theorem.
\end{pf}

Combining Theorem 7 and Theorem 8 we have the following:

\noindent{\it \textbf{Corollary 9.}} Let $\theta$ be any lacunary sequence. $\Delta N_{\theta} = |\Delta \sigma^{0}_{1}| $ if and only if $$1<lim\; inf \;q_{r}\leq lim\; sup \;q_{r} < \infty.$$

\noindent{\it \textbf{Corollary 10.}} Let $\theta$ be any lacunary sequence. $\Delta N_{\theta} = |\Delta \sigma^{0}_{1}|$ if and only if $|\sigma_{1}|=N_{\theta}$.

In the sequel, we will always assume that $\lim inf_{r}\; q_{r}>1$.

Now we give the definition of $N_{\theta}$-ward compactness.

\begin{Def}
 A subset $A$ of  $\Re$ is called $N_{\theta}$-ward compact if whenever $(\alpha _{n})$ is a sequence of points in $A$ there is a $N_{\theta}$-quasi-Cauchy subsequence $\textbf{z}=(z_{k})=(a_{n_{k}})$ of $(\alpha _{n})$.

\end{Def}

\noindent{\it \textbf{Theorem 11.}} A subset $A$ of $\Re$ is bounded if and only if it is $N_{\theta}$-ward compact.

\begin{pf}
Although there is a proof for the of lacunary statistical quasi-Cauchy sequences in an unpublished work of \c{C}akall\i, we give a proof for completeness.
Let $A$ be any bounded subset of $\Re$ and  $(\alpha_{n})$ be any sequence of points in $A$.  $(\alpha_{n})$  is also a sequence of points in $\overline{A}$  where  $\overline{A}$ denotes the closure of $A$. As  $\overline{A}$ is sequentially compact there is a convergent subsequence $(\alpha_{n_{k}})$ of $(\alpha_{n})$ (no matter the limit is in $A$ or not). This subsequence is $N_{\theta}$-convergent since $N_{\theta}$-method is regular. Hence $(\alpha_{n_{k}})$ is $N_{\theta}$-quasi-Cauchy. Thus (a) implies (b). To prove that (b) implies (a), suppose that $A$ is unbounded. If it is unbounded above, then one can construct a sequence $(\alpha _{n})$ of numbers in  $A$ such that  $\alpha_{n+1}>1+\alpha _{n}$ for each positive integer $n$. Then the sequence $(\alpha_{n})$ does not have any $N_{\theta}$-quasi-Cauchy subsequence, so $A$ is not $N_{\theta}$-ward compact. If $A$ is bounded above and unbounded below, then similarly we obtain that $A$ is not $N_{\theta}$-ward compact. This completes the proof.\end{pf}

\noindent{\it \textbf{Corollary 12.}} A closed subset of $\Re$ is $N_{\theta}$-ward compact if and only if it is $N_{\theta}$-sequentially compact.

\noindent{\it \textbf{Corollary 13.}} A subset of $\Re$ is $N_{\theta}$-ward compact if and only if it is statistically ward compact.

A sequence $\alpha=(\alpha _{n})$ is $\delta$-quasi-Cauchy if $\lim_{k\rightarrow \infty} \Delta^{2} \alpha _{n}=0$ where $\Delta^{2}\alpha _{n}=a_{n+2}-2a_{n+1}+\alpha _{n}$ (\cite{CakalliDeltaquasiCauchysequences}). A subset $A$ of $\Re$ is called $\delta$-ward compact if whenever $\boldsymbol \alpha=(\alpha _{n})$ is a sequence of points in $A$ there is a subsequence $\textbf{z}=(z_{k})=(\alpha_{n_{k}})$ of $\boldsymbol \alpha$ with $\lim_{k\rightarrow \infty} \Delta^{2}z_{k}=0$. It follows from the above theorem that any $N_{\theta}$-ward compact subset of $\Re$ is $\delta$-ward compact.

We note that $N_{\theta}$-quasi-Cauchy sequences were studied in \cite{BasarirandSelmaAltundagDeltaLacunarystatistical} in a different point of view.

We see that for any regular subsequential method $G$ defined on $\Re$, if a subset $A$ of $\Re$ is $G$-sequentially compact, then it is $N_{\theta}$- ward compact. But the converse is not always true.

\begin{Def}

 A subset $A$ of  $\Re$ is called strongly Cesaro compact if whenever $(\alpha _{n})$ is a sequence of points in $A$ there is a strongly Cesaro summable  subsequence $\textbf{z}=(z_{k})=(\alpha _{n_{k}})$ of $(\alpha _{n})$ to a strongly Cesaro limit that belongs to $A$.

\end{Def}

\begin{Def}

 A subset $A$ of  $\Re$ is called strongly Cesaro  ward compact if whenever $(\alpha _{n})$ is a sequence of points in $A$ there is a strongly Cesaro quasi-Cauchy subsequence $\textbf{z}=(z_{k})=(\alpha _{n_{k}})$ of $(\alpha _{n})$.

\end{Def}

\maketitle

\section{$N_{\theta}$-ward continuity}

\begin{Def}

A function defined on a subset $A$ of $\Re$ is called strongly Cesaro continuous if it preserves strongly Cesaro summable sequences, i.e. $(f(\alpha_{n}))$ is a strongly Cesaro quasi-Cauchy sequence whenever $(\alpha_{n})$ is.
\end{Def}

\begin{Def}

A function defined on a subset $A$ of $\Re$ is called strongly Cesaro ward continuous if it preserves strongly Cesaro quasi-Cauchy sequences, i.e. $(f(\alpha_{n}))$ is a strongly Cesaro quasi-Cauchy sequence whenever $(\alpha_{n})$ is.
\end{Def}

Now we give the definition of $N_{\theta}$-ward continuity in the following.

\begin{Def}
A function defined on a subset $A$ of $\Re$ is called $N_{\theta}$-ward continuous if it preserves $N_{\theta}$-quasi-Cauchy sequences, i.e. $(f(\alpha_{k}))$ is a $N_{\theta}$-quasi-Cauchy sequence whenever $(\alpha_{k})$ is.
\end{Def}

We note that $N_{\theta}$-ward continuity cannot be obtained by any sequential method $G$.
%(see \cite{CakalliOnGcontinuity} for more information on $G$-continuity).

A composition of two $N_{\theta}$-ward continuous functions is $N_{\theta}$-ward continuous. Sum of two $N_{\theta}$-ward continuous functions is $N_{\theta}$-ward continuous, but product of $N_{\theta}$-ward continuous functions need not be $N_{\theta}$-ward continuous.

In connection with $N_{\theta}$-quasi-Cauchy sequences and convergent sequences the problem arises to investigate the following types of  "continuity" of functions on $\Re$.

\begin{description}
\item[($\delta N_{\theta}$)] $(\alpha _{n}) \in {\Delta N^{0}_{\theta}} \Rightarrow (f(\alpha  _{n})) \in {\Delta N^{0}_{\theta}}$
\item[($\delta N_{\theta} c$)] $(\alpha  _{n}) \in {\Delta N^{0}_{\theta}} \Rightarrow (f(\alpha  _{n})) \in {c}$
\item[$(c)$] $(\alpha  _{n}) \in {c} \Rightarrow (f(\alpha  _{n})) \in {c}$
\item[$(c \delta N_{\theta})$] $(\alpha  _{n}) \in {c} \Rightarrow (f(\alpha  _{n})) \in {\Delta N^{0}_{\theta}}$
\item[($N_{\theta}$)] $(\alpha _{n}) \in {N^{0}_{\theta}} \Rightarrow (f(\alpha  _{n})) \in {N^{0}_{\theta}}$
\end{description}

We see that $(\delta N_{\theta})$ is $N_{\theta}$-ward continuity of $f$, $(N_{\theta})$ is $N_{\theta}$-continuity of $f$, and $(c)$ is the ordinary continuity of $f$. It is easy to see that $(\delta  N_{\theta} c)$ implies $(\delta N_{\theta})$, and $(\delta N_{\theta} )$ does not imply $(\delta N_{\theta} c)$;  and $(\delta N_{\theta})$ implies $(c \delta N_{\theta})$, and $(c \delta N_{\theta})$ does not imply $(\delta N_{\theta})$; $(\delta N_{\theta} c)$ implies $(c)$, and $(c)$ does not imply $(\delta N_{\theta} c)$; and $(c)$ is equivalent to $(c \delta N_{\theta})$.

Now we give the implication $(\delta N_{\theta})$ implies $(N_{\theta})$, i.e. any $N_{\theta}$-ward continuous function is $N_{\theta}$-continuous.

\noindent{\it \textbf{Theorem 14.}} If $f$ is $N_{\theta}$-ward continuous on a subset $A$ of $\Re$, then it is $N_{\theta}$-continuous on $A$.
\begin{pf}
Assume that $f$ is a $N_{\theta}$-ward continuous function on a subset $A$ of $\Re$.
Let $(\alpha _{n})$ be any $N_{\theta}$-convergent sequence with $N_{\theta}-\lim_{k\rightarrow\infty}\alpha _{k}=\alpha_{0}$. Then the sequence
$$(\alpha _{1}, \alpha _{0}, \alpha _{2}, \alpha _{0}, ..., \alpha _{n-1}, \alpha _{0}, \alpha _{n}, \alpha _{0},...)$$
is $N_{\theta}$-convergent to $\alpha _{0}$.  Hence it is $N_{\theta}$-quasi-Cauchy. As $f$ is $N_{\theta}$-ward continuous, the sequence  $$(f(\alpha _{1}),f(\alpha _{0}), f(\alpha _{2}), f(\alpha _{0}),..., f(\alpha _{n-1}), f(\alpha _{0}), f(\alpha _{n}), f(\alpha _{0}),...)$$ is $N_{\theta}$-quasi-Cauchy. It follows from this that the sequence $(f(\alpha _{n}))$ $N_{\theta}$-converges to $f(\alpha _{0})$. This completes the proof of the theorem.

The converse is not always true for the function $f(x)=x^{2}$ is an example since the sequence $(\sqrt{n})$ is $N_{\theta}$-quasi-Cauchy while $(f(\sqrt{n}))=(n)$ is not.
\end{pf}

It is well known that any continuous function on a compact subset $A$ of $\Re$ is uniformly continuous on $A$. It is also true for a regular subsequential method $G$ that any $N_{\theta}$-ward continuous function on a $G$-sequentially compact subset $A$ of $\Re$ is also uniformly continuous on $A$ (see [6]). Furthermore, for $N_{\theta}$-ward continuous functions defined on a $N_{\theta}$-ward compact subset of $\Re$, we have the following.

\noindent{\it \textbf{Theorem 15.}} Let $A$ be a $N_{\theta}$-ward compact subset $A$ of $\Re$ and let $f:A\longrightarrow$ $\Re$ be a  $N_{\theta}$-ward continuous function on $A$. Then $f$ is uniformly continuous on $A$.
\begin{pf}
Suppose that $f$ is not uniformly continuous on $A$ so that there exists an  $\varepsilon_{0} > 0$ such that for any $\delta >0$ $x, y \in{E}$ with $|x-y|<\delta$ but $|f(x)-f(y)| \geq \varepsilon_{0}$. For each positive integer $n$, fix $|\alpha _{n}-\beta_{n}|<\frac{1}{n}$, and $|f(\alpha _{n})-f(\beta_{n})|\geq \varepsilon_{0}$. Since $A$ is $N_{\theta}$-ward compact, there exists a $N_{\theta}$-quasi-Cauchy subsequence $(\alpha _{n_{k}})$ of the sequence $(\alpha _{n})$. It is clear that the corresponding subsequence $(\beta_{n_{k}})$ of the sequence $(\beta_{n})$ is also $N_{\theta}$-quasi-Cauchy, since $(\beta_{n_{k+1}}-\beta_{n_{k}})$ is a sum of three $N_{\theta}$-null sequences, i.e. $$\beta_{n_{k+1}}-\beta_{n_{k}}=(\beta_{n_{k+1}}-\alpha _{n_{k+1}})+(\alpha _{n_{k+1}}-\alpha _{n_{k}})+(\alpha _{n_{k}}-\beta_{n_{k}}).$$  On the other hand it follows from the equality $\alpha_{n_{k+1}}-\beta_{n_{k}}=\alpha_{n_{k+1}}-\alpha_{n_{k}}+\alpha_{n_{k}}-\beta_{n_{k}}$ that the sequence $(\alpha_{n_{k+1}}-\beta_{n_{k}})$ is $N_{\theta}$-convergent to $0$.
Hence the sequence
$$(a_{n_{1}}, \beta_{n_{1}}, \alpha _{n_{2}}, \beta_{n_{2}}, \alpha _{n_{3}}, \beta_{n_{3}},..., \alpha _{n_{k}}, \beta_{n_{k}},...)$$
is $N_{\theta}$-quasi-Cauchy.
But the transformed sequence
$$(f(\alpha _{n_{1}}), f(\beta_{n_{1}}), f(\alpha _{n_{2}}), f(\beta_{n_{2}}), f(\alpha _{n_{3}}), f(\beta_{n_{3}}),..., f(\alpha _{n_{k}}), f(\beta_{n_{k}}),...)$$
is not $N_{\theta}$-quasi-Cauchy. Thus $f$ does not preserve $N_{\theta}$-quasi-Cauchy sequences. This contradiction completes the proof of the theorem.
\end{pf}

\noindent{\it \textbf{Corollary 16.}} If a function $f$ is $N_{\theta}$-ward continuous on a bounded subset $A$ of $\Re$, then it is uniformly continuous on $A$.
\begin{pf}
The proof follows from the preceding theorem and Theorem 5.

\end{pf}

\noindent{\it \textbf{Theorem 17.}} $N_{\theta}$-ward continuous image of any $N_{\theta}$-ward compact subset of $\Re$ is $N_{\theta}$-ward compact.

\begin{pf}
Assume that $f$ is a $N_{\theta}$-ward continuous function on a subset $A$ of $\Re$, and $E$ is a $N_{\theta}$-ward compact subset of $A$. Let $(\beta _{n})$ be any sequence of points in $f(E)$. Write $\beta _{n}=f(\alpha _{n})$ where $\alpha _{n} \in {E}$ for each positive integer $n$. $N_{\theta}$-ward compactness of $E$ implies that there is a subsequence $(\gamma _{k})=(\alpha _{n_{k}})$ of $(\alpha _{n})$ with $N_{\theta}\lim_{k\rightarrow \infty} \Delta \gamma _{k}=0$. Write $(t_{k})=(f(\gamma_{k}))$. As $f$ is $N_{\theta}$-ward continuous, $(f(\gamma_{k}))$ is $N_{\theta}$-quasi-Cauchy. Thus we have obtained a subsequence $(t_{k})$ of the sequence $(f(\alpha _{n}))$ with $N_{\theta}-\lim_{k\rightarrow \infty} \Delta t_{k}=0$. Thus $f(E)$ is $N_{\theta}$-ward compact. This completes the proof of the theorem.
\end{pf}

\noindent{\it \textbf{Corollary 18.}} $N_{\theta}$-ward continuous image of any compact subset of $\Re$ is $N_{\theta}$-ward compact.

The proof follows from the preceding theorem.

\noindent{\it \textbf{Corollary 19.}} $N_{\theta}$-ward continuous image of any bounded subset of $\Re$ is bounded.

The proof follows from Theorem 6 and Theorem 15.

\noindent{\it \textbf{Corollary 20.}} $N_{\theta}$-ward continuous image of a $G$-sequentially compact subset of $\Re$ is $N_{\theta}$-ward compact for any
  regular subsequential method $G$.

It is a well known result that uniform limit of a sequence of continuous functions is continuous. This is also true in case of $N_{\theta}$-ward continuity, i.e. uniform limit of a sequence of $N_{\theta}$-ward continuous functions is $N_{\theta}$-ward continuous.

\noindent{\it \textbf{Theorem 21.}} If $(f_{n})$ is a sequence of $N_{\theta}$-ward continuous functions on a subset $A$ of $\Re$ and $(f_{n})$ is uniformly convergent to a function $f$, then $f$ is $N_{\theta}$-ward continuous on $A$.
\begin{pf}  Let $\varepsilon$ be a positive real number and $(\alpha _{k})$ be any $N_{\theta}$-quasi-Cauchy sequence of points in $A$. By uniform convergence of $(f_{n})$ there exists a positive integer $N$ such that $|f_{k}(x)-f(x)|<\frac{\varepsilon}{4}$ for all $x \in {E}$ whenever $k\geq N$. As $f_{N}$ is $N_{\theta}$-ward continuous on $A$, there exists an $n_{0}$, greater than $N$, such that $r>n_{0}$ implies that
$\frac{1}{h_{r}}\sum^{}_{k\in{I_{r}}}|\Delta f_{N}(\alpha_{k}-\alpha_{k+1}|<\frac{\varepsilon}{2}$. Thus for $k>n$ we have \\
$\frac{1}{h_{r}}\sum^{}_{k\in{I_{r}}}|\Delta f(\alpha_{k})| \leq {\frac{1}{h_{r}}\sum^{}_{k\in{I_{r}}}|f(\alpha_{k})-f_{N}(\alpha_{k})|+\frac{1}{h_{r}}\sum^{}_{k\in{I_{r}}}|f_{N}(\alpha_{k+1})- f(\alpha_{k+1})|+\frac{1}{h_{r}}\sum^{}_{k\in{I_{r}}}|f_{N} (\alpha_{k})-f_{N}(\alpha_{k+1})|}< \frac{1}{h_{r}} (k_{r}-k_{r-1})\frac{\varepsilon}{4}+\frac{\varepsilon}{2}=\varepsilon$\\
This completes the proof of the theorem.

\end{pf}

\noindent{\it \textbf{Theorem 22.}} The set of all $N_{\theta}$-ward continuous functions on a subset $A$ of $\Re$ is a closed subset of the set of all continuous functions on $A$, i.e. $\overline{\Delta N_{\theta}WC(A)}=\Delta N_{\theta}WC(A)$ where $\Delta N_{\theta}WC(A)$ is the set of all $N_{\theta}$-ward continuous functions on $A$, $\overline{\Delta N_{\theta}WC(A)}$ denotes the set of all cluster points of $\Delta N_{\theta}WC(A)$.
\begin{pf}
The proof of this theorem is sþmilar to that of the preceding theorem so is omitted.
%Let $f$ be any element in the closure of $\Delta N_{\theta}WC(A)$. Then there exists a sequence of points in $\Delta N_{\theta}WC(A)$ such that $\lim_{k\rightarrow \infty} f_{k}=f$. To show that $f$ is $N_{\theta}$-ward continuous take any $N_{\theta}$-quasi-Cauchy sequence $(\alpha _{k})$ of points in $A$. Let $\varepsilon > 0$. Since $(f_{k})$ converges to $f$, there exists an $N$ such that for all $x \in {A}$ and for all $n \geq {N}$, $|f(x)-f_{n}(x)|< \frac{\varepsilon}{3}$. Since $f_N$ is $N_{\theta}$-ward continuous on $A$, there exists an $n_{0}$, greater than $N$, such that $r>n_{0}$ implies that $\frac{1}{h_{r}}\sum^{}_{k\in{I_{r}}}|\Delta f_{N}(\alpha_{k}-\alpha_{k+1}|<\frac{\varepsilon}{2}$. Thus for $k>n$ we have \\ $\frac{1}{h_{r}}\sum^{}_{k\in{I_{r}}}|\Delta f(\alpha_{k})| \leq {\frac{1}{h_{r}}\sum^{}_{k\in{I_{r}}}|f(\alpha_{k})-f_{N}(\alpha_{k})|+\frac{1}{h_{r}}\sum^{}_{k\in{I_{r}}}|f_{N}(\alpha_{k+1})- f(\alpha_{k+1})|+\frac{1}{h_{r}}\sum^{}_{k\in{I_{r}}}|f_{N} (\alpha_{k})-f_{N}(\alpha_{k+1})|}< \frac{1}{h_{r}} (k_{r}-k_{r-1})\frac{\varepsilon}{4}+\frac{\varepsilon}{2}=\varepsilon$\\ This completes the proof of the theorem.
\end{pf}
\noindent{\it \textbf{Corollary 23.}} The set of all $N_{\theta}$-ward continuous functions on a subset $A$ of $\Re$ is a complete subspace of the space of all continuous functions on $A$.
\begin{pf}
The proof follows from the preceding theorem.
\end{pf}

\maketitle

\section{Conclusion}
In this paper, new types of continuities are introduced via $N_{\theta}$-quasi-Cauchy sequences and strongly Cesaro quasi-Cauchy sequences, and investigated. In the investigation we have obtained results related to $N_{\theta}$-ward continuity, $N_{\theta}$-continuity, strongly Cesaro ward continuity, strongly Cesaro continuity, $G$-sequential continuity, ordinary continuity, uniform continuity, $N_{\theta}$-ward compactness, $N_{\theta}$-sequentially compactness, statistical ward compactness, statistical compactness, $G$-sequential compactness, and ordinary compactness. We also proved some inclusion theorems between the set of $N_{\theta}$-quasi-Cauchy sequences and the set of statistical quasi-Cauchy sequences.

For further study, we suggest to investigate $N_{\theta}$-quasi-Cauchy sequences of fuzzy points, strongly Cesaro quasi-Cauchy sequences, and $N_{\theta}$-ward continuity for the fuzzy functions (see \cite{CakalliandPratul} for the definitions and  related concepts in fuzzy setting). However due to the change in settings, the definitions and methods of proofs will not always be analogous to those of the present work.

\end{document}